\documentclass[11pt,a4paper]{article}

\usepackage{amsmath, amsfonts, amssymb}
\usepackage{color}

\def\be#1\ee{\begin{equation}#1\end{equation}}

\newtheorem{thm}{Theorem}
\newtheorem{lem}[thm]{Lemma}

\def\P{{\mathbb{P}}}
\def\R{\mathbb{R}}
\def\E{\mathbb{E}\,}

\def\N{{\mathbb N}}
\def\dd{\mbox{d}}

\def\d{\, \mathrm{d}}

\newcommand{\ind}{1\hspace{-0.098cm}\mathrm{l}}

\newcommand{\eps}{\varepsilon}

\let\BFseries\bfseries\def\bfseries{\BFseries\mathversion{bold}} 

\def\WDHH{\mathbf{W}^H_{d+1}}
\def\WDH{\mathbf{W}^H_d}
\def\bfX{\mathbf{X}}
\def\bfY{\mathbf{Y}}
\def\bfR{\mathbf{R}}
\def\WH{W^H}
\def\WTH{W^{\tilde H}}
\def\F{\mathcal{F}}
\def\deq{{~\stackrel{d}{=}~}}

\title{The first exit time of fractional Brownian motion from a parabolic domain}
\author{Frank Aurzada and Mikhail Lifshits}

\begin{document}

 \maketitle


\begin{abstract}
 We study the first exit time of a multi-dimensional fractional Brownian motion from
 unbounded domains. In particular, we are interested in the upper tail of
 the corresponding distribution when the domain is parabola-shaped.
\end{abstract}



\section{Introduction and main result}
Let $\WDHH$ be a $(d+1)$-dimensional fractional Brownian motion (FBM), i.e a process with independent
coordinates, each representing a classical one-dimensional FBM with parameter $0<H<1$. We are interested in the
first exit time of this process from an unbounded domain $D$:
$$
    \tau_D:=\inf\{ t > 0 : \WDHH(t) \not\in D \}.
$$
To the knowledge of the authors, very little is known even for very specific domains
$D$ like cones, etc, except for the Brownian motion case ($H=1/2$).

A specially interesting setting appears when $D$ is a cone; here we refer to
\cite{spitzer, deblassie87, banuelos-smits,banuelos-deblassie-smits} for results in the Brownian case. The case of
multi-dimensional fractional Brownian motion (even two-dimensional) has not been considered
and is an interesting open question. We comment on the case of cones in Section~\ref{sec:generalizations} below.

 In this paper, we focuss on parabola shaped domains
$$
D:= \{ (x_1,\ldots,x_{d+1}) : || (x_1,\ldots,x_d) ||^p \leq a + x_{d+1} \}
$$
with parameters $p>0$ and $a>0$. Here and throughout, $||.||$ denotes the standard Euclidian norm.
This paper is a companion to a work by Lifshits and Shi \cite{lifshitsshi} where the Brownian
case is treated  (also see \cite{li2003} and \cite{banuelos-deblassie-smits} for earlier results).
The mentioned paper \cite{lifshitsshi} was used and extended by many authors to various
processes, see e.g. \cite{nane} for iterated Brownian motion and \cite{banuelosbogdan,banuelosbogdan2} for
symmetric stable processes.
However, to our knowledge the case of multi-dimensional fractional Brownian motion has not been
considered so far.

We shall see that the tail probabilities of $\tau_D$  are governed by a mixture of small and
large deviations: the most favorable way for $\tau_D$ to be large is when the first $d$
components are in a certain small deviation regime, while the $(d+1)$-th component performs
in a large deviation mode. In order to state the result, we need to introduce the small deviation
constant for multivariate FBM (see e.g.\ \cite{LiLinde}):
\begin{equation}
\label{eqn:sdconstantfbm}
     \lim_{\eps\to\infty} - \eps^{1/H} \log \P( \sup_{t\in[0,1]} ||\WDH(t)||\leq \eps )
     =: \kappa_{H,d} \in(0,\infty).
\end{equation}

We are now able to formulate the main result of this paper.

\begin{thm} \label{thm:mainthm}
Let $p>1$ and set $\beta:=2H(p-1)/(2Hp+1)$. Let $\WDH$ be a $d$-dimensional fractional Brownian motion
and $\WH$ a one-dimensional fractional Brownian motion, independent of $\WDH$, let $a>0$. Then
$$
    \lim_{T\to\infty} T^{-\beta} \log \P( ||\WDH(t) ||^p \leq a + \WH(t), \forall t\in[0,T] )
    = -\kappa \in (-\infty,0),
$$
where $\kappa$ is the solution of the following minimization problem
$$
    \kappa:=\inf_{h(\cdot) \geq 0} \left[ \kappa_{H,d} \int_0^1 h(t)^{-1/(pH)}\d t + \frac{||h||_H^2}{2}\right],
$$
and $||.||_H$ is the RKHS norm related to FBM with parameter $H$.
\end{thm}

In the case $p<1$, the probability in question decays as $T^{-(1-H)+o(1)}$; we comment on this fact 
in Section~\ref{sec:generalizations}. In the case $p=1$, the domain is in fact a cone. Using our methods, 
one can show that the probability also decays polynomially, but we are not able to determine the exact 
exponent, see Section~\ref{sec:generalizations}.

The rest of the paper is arranged as follows: In Section~\ref{sec:lower}, we give the proof of the lower 
bound in Theorem~\ref{thm:mainthm}, while the proof of the corresponding upper bound is given in Section~\ref{sec:upper}. 
The final Section~\ref{sec:generalizations} contains some comments on possible generalizations and open questions.

\section{Lower bound} \label{sec:lower}

Let $h(t):=T^q h_0(t/T)$, where $q:=(2H+1)Hp/(2Hp+1)$ and $h_0$ is some fixed non-negative
function from the RKHS of FBM. It is simple to check that $p>1$ implies $q\in(H,pH)$.
We have
\begin{eqnarray}
&& \P( ||\WDH(t) ||^p \leq a + \WH(t), \forall t\in[0,T] )
\notag
\\
&\geq & \P( ||\WDH(t) ||^p \leq h(t)\leq a + \WH(t), \forall t\in[0,T] )
\notag
\\
&= & \P( ||\WDH(t) ||^p \leq h(t), \forall t\in[0,T]) \cdot \P(h(t) \leq a + \WH(t), \forall t\in[0,T] )
\label{eqn:lowerboundsplit}
\end{eqnarray}

The first probability is a small deviation probability, the second will turn
out to be a large deviation probability. Let us start with the first probability:
Using the self-similarity of $\WDH$, we have
\begin{eqnarray*}
&& \P( ||\WDH(t) ||^p \leq h(t), \forall t\in[0,T])
\\
&=&
\P( ||T^H \WDH(t/T) || \leq T^{q/p} h_0(t/T)^{1/p}, \forall t\in[0,T])
\\
&=&
\P( ||\WDH(t) || \leq T^{q/p-H} h_0(t)^{1/p}, \forall t\in[0,1]).
\end{eqnarray*}

Since $q/p-H<0$, we can apply a formula from \cite{lifshitslindetams} (cf.\ Theorem~4.6 and the remark on the uniform norm on p.~2072 there;
notice that the results of \cite{lifshitslindetams} handle one-dimensional processes but the multivariate ones easily follow by the same
methods based on subadditivity arguments of \cite{LiLinde}) and see that this probability
admits a representation
$$
\exp\left( - \kappa_{H,d} \int_0^1 h_0(t)^{-1/(pH)} \d t \cdot T^{-\frac{1}{H}(q/p-H)} (1+o(1))\right).
$$
where $\kappa_{H,d}$ is as in (\ref{eqn:sdconstantfbm}). A quick computation shows that $-\frac{1}{H}(q/p-H)=\beta$.

We now deal with the second probability in (\ref{eqn:lowerboundsplit}): using the self-similarity of $\WH$, we get
\begin{eqnarray*}
&& \P(h(t) \leq a + \WH(t), \forall t\in[0,T] )
\\
&=& \P(T^q h_0(t/T) \leq a + T^H \WH(t/T), \forall t\in[0,T] )
\\
&=& \P( \WH(t) \geq - a T^{-H} + T^{q-H} h_0(t) , \forall t\in[0,1] ).
\end{eqnarray*}

At this point, we apply Proposition~1.6 from \cite{aurzadadereich}: It states that for any $f_0$ from the reproducing kernel
Hilbert space of a Gaussian random vector $X$ attaining values in some Banach space and for any measurable set $S$ one has
$$
    \P( X + f_0 \in S) \geq \P( X \in S)\cdot
    \exp\left( - \frac{||f_0||_{\mathcal{H}}^2}{2} - \sqrt{2 ||f_0||_{\mathcal{H}}^2 \log \P( X\in S)^{-1}} \right),
$$
where $||.||_{\mathcal{H}}$ is the norm in the reproducing kernel Hilbert space. A similar upper bound holds,
but we do not need it in our context. We use this inequality with the set $S:=\{ f : f(t)\geq -a T^{-H}, t\in[0,1]\}$
 and the function $f_0=-T^{q-H} h_0$, which clearly belongs to the RKHS of $\WH$. The inequality then gives
\begin{eqnarray*}
      && \P( \WH(t) \geq -a T^{-H} + T^{q-H} h_0(t) , \forall t\in[0,1] )
\\
     &\geq& p_T\cdot \exp\left( -|| T^{q-H} h_0||_H^2/2 - \sqrt{2 || T^{q-H} h_0||_H^2 \log p_T^{-1}} \right),
\end{eqnarray*}
where $p_T=\P( \WH(t) \geq - a T^{-H}, \forall t\in[0,1] )=T^{-(1-H)+o(1)}$, see e.g.\ \cite{molchan}. This shows that
\begin{eqnarray*}
    \P( \WH(t) &\geq& -a T^{-H} + T^{q-H} h_0(t) , \forall t\in[0,1] )
\\
    &=& \exp\left( - T^{2q-2H} \frac{||h_0||_H^2}{2} \cdot(1+o(1)) \right).
\end{eqnarray*}
Another quick computation yields that also $2q-2H=\beta$.

Putting everything together, we obtain
\begin{eqnarray}
&& \liminf_{T\to\infty} T^{-\beta} \log \P( ||\WDH(t) ||^p \leq a + \WH(t), \forall t\in[0,T] )
\notag
\\
&\geq&
- \left( \kappa_{H,d} \int_0^1 h_0(t)^{-1/(pH)} \d t + \frac{||h_0||_H^2}{2}\right), \label{eqn:lowerbnot-inf}
\end{eqnarray}
for any non-negative function $h_0$ from the RKHS of FBM. Taking the supremum over all such functions
shows the lower bound in Theorem~\ref{thm:mainthm}.

\section{Upper bound} \label{sec:upper}

We are now able to prove the upper bound in Theorem \ref{thm:mainthm}.
Our approach here is based on the classical Varadhan lemma from large deviation theory,
see Step~4 in the proof below.


{\it Step 1: Discretization.}

Let $0=t_0< t_1 < \ldots < t_N = T$ with $N\in\N$ and observe that
\begin{eqnarray}
&& \P( ||\WDH(t) ||^p \leq a + \WH(t), \forall t\in[0,T] ) \notag
\\
   &\leq& \P( \sup_{t\in[t_{k-1},t_k]} ||\WDH(t) || \leq (a + \sup_{t\in[t_{k-1},t_k]} \WH(t))_+^{1/p}, \forall k\leq N ) \notag
\\
   &=&\E_{\WH} \P( \sup_{t\in[t_{k-1},t_k]} ||\WDH(t) || \leq a_k, \forall k\leq N ),
\label{eqn:cutofftrick}
\end{eqnarray}
where $a_k:=(a + \sup_{t\in[t_{k-1},t_k]} \WH(t))_+^{1/p}$.

\noindent{\it Step 2: Evaluation of a small deviation probability.}

Now recall that for $t\in[t_{N-1},t_N]$,  $\WDH(t)=(W_H^{(1)}(t),\ldots,W_H^{(d)}(t))$
can be rewritten, by the Mandelbrot-van Ness representation, as follows:
\begin{eqnarray*}
 W_H^{(i)}(t)  &=&\int_{-\infty}^{t_{N-1}} \left[ (t-u)^{H-1/2} -(-u)_+^{H-1/2} \right]  \dd W^{(i)}(u)
\\
   &&    + \int_{t_{N-1}}^t (t-u)^{H-1/2}  \dd W^{(i)}(u),
\\
   &=:& X^{(i)}(t) + Y^{(i)}(t),
\end{eqnarray*}
where the processes $W^{(i)}$ are independent Brownian motions.
Writing $\bfX=(X^{(1)},\ldots,X^{(d)})$ and similarly for $\bfY$, note that the probability
on the right hand side in (\ref{eqn:cutofftrick}) can be rewritten as
$$
   \E[ \ind_{\{\sup_{t\in[t_{k-1},t_k]} ||\WDH(t) || \leq a_k, \forall k\leq N-1\}}
   \P( \sup_{t\in[t_{N-1},t_N]} ||\bfX(t)+\bfY(t) || \leq a_N | \F_{N-1})],
$$
where $\F_{N-1}:=\sigma(W^{(i)}(t), t\leq t_{N-1}, i=1,\ldots,d)$. Noting that $\bfX$ on $[t_{N-1},t_N]$
is determined by $\F_{N-1}$, one can eliminate it by Anderson's inequality,
so that the last term can be estimated from above by
$$
   \E[ \ind_{\{\sup_{t\in[t_{k-1},t_k]} ||\WDH(t) || \leq a_k, \forall k\leq N-1\}}
   \P( \sup_{t\in[t_{N-1},t_N]} ||\bfY(t) || \leq a_N | \F_{N-1})].
$$
We notice that $\bfY$ on $[t_{N-1},t_N]$ does not depend on $\F_{N-1}$,
so that the expression becomes
$$
   \P( \sup_{t\in[t_{k-1},t_k]} ||\WDH(t) || \leq a_k, \forall k\leq N-1)
   \cdot \P( \sup_{t\in[t_{N-1},t_N]} ||\bfY(t) || \leq a_N ).
$$
We further note that
\begin{eqnarray*}
    Y^{(i)}(t)
    &=& \int_{t_{N-1}}^t (t-u)^{H-1/2}  \dd W^{(i)}(u)
\\
   &\deq& \int_{0}^{t-t_{N-1}} (t-t_{N-1}-u)^{H-1/2}  \dd W^{(i)}(u)
\\
   &=:& R_N^{(i)}(t-t_{N-1}),
\end{eqnarray*}
where $\deq$ is means the equality of finite dimensional distributions.
The processes $R^{(i)}_N$ are called Riemann-Liouville processes. We
write $\bfR_N:=(R^{(1)}_N,\ldots,R^{(d)}_N)$ for short.

Iterating the above arguments, one obtains
$$
    \P( \sup_{t\in[t_{k-1},t_k]} ||\WDH(t) || \leq a_k, \forall k\leq N )
    \leq \prod_{k=1}^N \P( \sup_{t\in[t_{k-1},t_k]} ||\bfR_k(t-t_{k-1}) || \leq a_k ),
$$
where the $\bfR_k$ are independent copies of $\bfR_N$.

By the self-similarity of the $\bfR_k$, the last expression is equal to
\begin{equation} \label{eqn:nowsmalldevations}
    \prod_{k=1}^N \P( \sup_{t\in[0,t_k-t_{k-1}]} ||\bfR_k(t) || \leq a_k )
    = \prod_{k=1}^N \P( \sup_{t\in[0,1]} ||\bfR_k(t) || \leq \frac{a_k}{(t_k-t_{k-1})^H} ).
\end{equation}

Now we recall the small deviation bounds for Riemann-Liouville processes from \cite{LiLinde} (also see \cite{ls}):
For any $\delta>0$ there is a constant $c>0$ such that for any $r>0$

$$
  \P( \sup_{t\in[0,1]} ||\bfR_{k}(t)|| \leq r ) \leq c \exp\left(- \kappa_{H,d}(1-\delta)r^{-1/H}  \right),
$$
where $\kappa_{H,d}$ is the same constant as for FBM in (\ref{eqn:sdconstantfbm}).

This implies that the product in (\ref{eqn:nowsmalldevations}) can be majorated by
$$
     c^N \exp\left( - \kappa_{H,d}(1-\delta) \sum_{k=1}^N \frac{t_k-t_{k-1}}{a_k^{1/H}} \right).
$$

Putting this together with (\ref{eqn:cutofftrick}), we obtain
\begin{eqnarray*}
   && \P( ||\WDH(t) ||^p \leq 1 + \WH(t), \forall t\in[0,T] )
\\
&\leq&
    c^N \E \exp\left( - \kappa_{H,d}(1-\delta) \sum_{k=1}^N
    \frac{t_k-t_{k-1}}{(1 + \sup_{t\in[t_{k-1},t_k]} \WH(t))_+^{1/(pH)}} \right).
\end{eqnarray*}

\noindent{\it Step 3: Simplification of the functional.}

Setting $\tau_k:= t_k/T$ and using that $(T^H \WH(t/T))$ has the same
finite dimensional distributions as $\WH$, the last inequality becomes
\begin{eqnarray*}
   && \P( ||\WDH(t) ||^p \leq a + \WH(t), \forall t\in[0,T] )
\\
   &\leq& c^N \E \exp\left( - \kappa_{H,d}(1-\delta) T^{1-1/p}
        \sum_{k=1}^N \frac{\tau_k-\tau_{k-1}}{(aT^{-H} + \sup_{t\in[\tau_{k-1},\tau_k]} \WH(t))_+^{1/(pH)}} \right).
\end{eqnarray*}

We now show how to get rid of the $aT^{-H}$ term. Fix $\theta>0$ and define the event
$E:=\{ \forall k =1,\ldots, N :  T^{-H} \leq \theta \sup_{t\in[\tau_{k-1},\tau_k]} \WH(t)\}$.
On the complementary event $E^c$, we have
\begin{eqnarray*}
    && \sum_{k=1}^N\frac{\tau_k-\tau_{k-1}}{(aT^{-H} + \sup_{t\in[\tau_{k-1},\tau_k]} \WH(t))_+^{1/(pH)}}
\\
    &\geq& \frac{\min_k(\tau_{k}-\tau_{k-1})}{T^{-1/p} (a + 1/\theta)^{1/(pH)}}
    =: \frac{m_\tau}{T^{-1/p} (a + 1/\theta)^{1/(pH)}}.
\end{eqnarray*}

Therefore, using the notation
\[
   S[f,\tau]:= \sum_{k=1}^N
           \frac{\tau_k-\tau_{k-1}}{(\sup_{t\in[\tau_{k-1},\tau_k]} f(t))_+^{1/(pH)}}
\]
for a function $f$ and a partition $\tau$, we obtain
\begin{eqnarray*}
     && \E \exp\left( - \kappa_{H,d}(1-\delta) T^{1-1/p} \sum_{k=1}^N \frac{\tau_k-\tau_{k-1}}
        {(aT^{-H} + \sup_{t\in[\tau_{k-1},\tau_k]} \WH(t))_+^{1/(pH)}} \right)
\\
    &\leq& \E \left[ \ind_E \exp\left( - \kappa_{H,d}\frac{1-\delta}{(a\theta+1)^{1/(pH)}} T^{1-1/p}
    S[\WH,\tau]  \right) \right]
\\
    &&  + \E \left[ \ind_{E^c} \exp\left( - \kappa_{H,d}(1-\delta) T \frac{m_\tau}{(a + 1/\theta)_+^{1/(pH)}} \right) \right]
\\
    &\leq& \E \exp\left( - \kappa_{H,d}\frac{1-\delta}{(1 + a\theta)^{1/(pH)}} T^{1-1/p}
    S[\WH,\tau]  \right)  + \exp\left( - c' T \right),
\end{eqnarray*}
for some non-random constant $c'$ depending only on the constants $\kappa_{H,d},\delta,\theta,a$ and on the
partition $(\tau_k)$ but not depending on $T$. We shall see below that the second term is of lower order,
as when taking logarithms, the order of the first term is $T^\beta$, with $\beta<1$.

\noindent{\it Step 4: Application of Varadhan's lemma.}

Recall a special case of Varadhan's lemma (see e.g.\ Theorem III.13 in \cite{denhollander} or Theorem 4.3.1 in \cite{dembozeitouni}):

\begin{lem}[Varadhan's lemma] Let $(Z_\eps)$ be a familiy of random elements of $\mathcal{C}[0,1]$.
Assume that the family $(\P_{Z_\eps})$ satisfies a large deviation principle with good rate funtion $\mathcal{I}$.
Let $\phi : \mathcal{C}[0,1] \to \R$ be a continuous function that is bounded from above. 
Then
\begin{equation}\label{eqn:statemevl}
   \lim_{\eps\to 0} \eps \log \E \exp\left(  \phi( Z_\eps)/\eps  \right) = \sup_f (\phi(f) -I(f)).
\end{equation}
\end{lem}

Note that the expression
\begin{equation} \label{eqn:prefac}
    \E \exp\left( - \kappa_{H,d}\frac{1-\delta}{(1 + a\theta)^{1/(pH)}} T^{1-1/p}
       S[\WH,\tau]\right)
\end{equation}
fits into the framework of Varadhan's lemma. The functional $\phi$ is given by
$$
     \phi(f):=- \kappa_{H,d}\frac{1-\delta}{(1 + a\theta)^{1/(pH)}} S[f,\tau]
$$
while we choose $\eps=T^{-\beta}$. Then the term in (\ref{eqn:prefac}) becomes
$$
    \E[ e^{\phi( \eps^{1/2} \WH )/\eps} ],
$$
and we know from the large deviation theory for Gaussian processes, see e.g.\ \cite[Chapter 12]{Lif}, that $Z_\eps=\eps^{1/2} \WH$ satisfies the large deviation principle
with good rate function $\mathcal{I}(f) := ||f||_H^2/2$, where $||.||_H$ is the RKHS norm.

 This means that (\ref{eqn:statemevl})
holds for FBM $Z_\eps=(\eps^{1/2} \WH(t))_{t\in[0,1]}$ and the continuous, bounded functional $\phi$
and we deduce that
\begin{eqnarray}
    && \limsup_{T\to\infty} T^{-\beta} \log \P( ||\WDH(t) ||^p \leq 1 + \WH(t), \forall t\in[0,T] )
    \notag
\\
   &\leq& \lim_{T\to\infty} T^{-\beta} \log \E \exp\left( - \kappa_{H,d}\frac{1-\delta}{(1 + a\theta)^{1/(pH)}}
   T^{1-1/p}  S[\WH,\tau] \right)
   \notag
\\
   &=& \sup_f (\phi(f) -\mathcal{I}(f))
   \notag
\\
   &=& - \inf_{f} \left[ \kappa_{H,d}\frac{1-\delta}{(1 + a\theta)^{1/(pH)}}  S[f,\tau] +||f||_H^2/2\right]
   =: - V_{\tau,\delta,\theta}.
   \label{eqn:defV}
\end{eqnarray}

\noindent{\it Step 4: Final computations.} Our final goal is to get rid of the partitions.

For this purpose, let $\rho_n\to 0$, $\delta_n\to 0$, $\theta_n\to 0$ and let $(\tau^{(n)})$
be a sequence of partitions with diameter tending to zero. Certainly, for each $n$ one can find
a function $f_n$ such that
$$
   \kappa_{H,d} \frac{1-\delta_n}{(1 + a\theta_n)^{1/(pH)}}
    S[f_n,\tau^{(n)}]    +||f_n||_H^2/2
   \leq V_{(\tau^{(n)}_k),\delta_n,\theta_n} + \rho_n.
$$

By the lower bound (\ref{eqn:lowerbnot-inf}), the sequence of functions must be bounded in $||.||_H$.
By the compactness of the RKHS balls, $(f_n)$ has a convergent subsequence (w.l.o.g.\ the original sequence);
let us denote the limit by $f$.
Further, since the diameters of the partitions tend to zero, we have
$$
    S[f_n,\tau^{(n)}]  \to \int_0^1 f(t)_+^{-1/(pH)}\d t.
$$

Therefore,
\begin{eqnarray*}
   V_{\tau^{(n)},\delta_n,\theta_n}
   &\geq& \kappa_{H,d}\frac{1-\delta_n}{(1 + a\theta_n)^{1/(pH)}}
    S[f_n,\tau^{(n)}]  +||f_n||_H^2/2 - \rho_n
\end{eqnarray*}
When letting $n\to\infty$, the bound tends to $\kappa_{H,d} \int_0^1 f(t)_+^{-1/(pH)}\d t + ||f||_H^2/2$,
which is, in turn, minorated by
$$
 \inf_{h \geq 0}  \left[ \kappa_{H,d} \int_0^1 h(t)^{-1/(pH)} \d t+ ||h||_H^2/2\right],
$$
which after putting it together with \eqref{eqn:defV} confirms the upper bound in Theorem~\ref{thm:mainthm}.

\section{Generalizations and open questions} \label{sec:generalizations}
\paragraph*{Distinct Hurst parameters.} Let us consider the following generalization of our problem. Let $\WDH$ be, as before, a $d$-dimensional FBM
with Hurst paramter $H$. Let $\WTH$ be a one-dimensional FBM with Hurst parameter $\tilde H$, independent
of $\WDH$. Let $p>0$. Then the probability in question
$$
\P( ||\WDH(t)||^p \leq a + \WTH(t), \forall t\in[0,T])
$$
admits three regimes: $pH>\tilde H$, $pH=\tilde H$, and $pH<\tilde H$.

If $p H>\tilde H$, everything works as in Theorem~\ref{thm:mainthm} and one obtains
$$
\lim_{T\to\infty} T^{-\tilde\beta} \log \P( ||\WDH(t) ||^p \leq a + \WH(t), \forall t\in[0,T] )
    = -\tilde\kappa \in (-\infty,0)
$$
with $\tilde\beta:=2(Hp-\tilde H)/(2Hp+1)$ and $\tilde\kappa$ is again given by a similar minimization problem.

If $\tilde H > pH$, the rate of decay becomes polynomial:
$$
\P( ||\WDH(t) ||^p \leq a + \WH(t), \forall t\in[0,T] ) = T^{-(1-\tilde H)+o(1)}.
$$

The proofs of these facts go along the lines of the proof of Theorem~\ref{thm:mainthm}
(choose $h(t)=T^q h_0(t/T)$ with $h_0$ in the reproducing kernel Hilbert space of $\WTH$ and
$q\in(\tilde H,pH)$ in the first case and $q\in(pH,\tilde H)$ in the second case).

The third regime is the critical case $\tilde H=pH$, which deserves more comments.

\paragraph*{Open critical cases.}
The critical case that remains open is  $\tilde H=pH$.
We believe that in this case the problem should be stated using a width parameter
$K$, namely, we conjecture
$$
    \P( ||\WDH(t) ||^p \leq K(a + \WTH(t)), \forall t\in[0,T] )
   = T^{-\gamma(K) +o(1)},
$$
i.e.\ the tail probability decays polynomially with a power depending on the width
(unlike in the subcritical case mentioned above). Our methods allow to show that the rate lies between polynomial functions, but we are not able to determine $\gamma(K)$.

In particular, when $\tilde H=H$ and $p=1$ the event in question actually concerns staying in
a {\it cone}. For Brownian motion, the corresponding tail probabilities were studied by Spitzer \cite{spitzer}
who obtained for $d=2$
\[
   \gamma(K)= \frac{\pi}{4\arctan(K)},
\]
and by DeBlassie   \cite{deblassie87,deblassie88} for general $d$ with an inexplicit
representation of $\gamma(K)$; see also  Ba\~nuelos and Smits \cite{banuelos-smits}.

\paragraph*{Width parameter.}
Also in the case studied in Theorem~\ref{thm:mainthm}, one can introduce a width parameter $K>0$.
A straightforward argument using the scaling properties of $\WDH$ and $\WH$ gives
\begin{eqnarray*}
    &&
    \lim_{T\to\infty} T^{-\beta} \log \P( ||\WDH(t) ||^p \leq K(a + \WH(t)), \forall t\in[0,T] )
    \\
    &=& -\kappa K^{-2/(2Hp+1)} \in (-\infty,0),
\end{eqnarray*}
where $\beta$ and $\kappa$ are the same as in Theorem~\ref{thm:mainthm}. The same assertion
holds for the generalized problem:
\begin{eqnarray*}
    && \lim_{T\to\infty} T^{-\tilde\beta} \log \P( ||\WDH(t) ||^p \leq K(1 + \WTH(t)), \forall t\in[0,T] )
    \\
    &=& -\tilde\kappa K^{-2/(2Hp+1)} \in (-\infty,0),
\end{eqnarray*}
with $\tilde\beta:=2(Hp-\tilde H)/(2Hp+1)$ and $\tilde\kappa$ as above.

\bigskip
{\bf Acknowledgement.} This research was supported by the Russian Foundation
Basic Research grant 16-01-00258 and by the co-ordinated
grants of DFG (GO420/6-1) and St.\ Petersburg State University (6.65.37.2017).

%
%

\noindent {\bf Addresses of the authors:}\\
Frank Aurzada, Technische Universit\"at Darmstadt, Schlossgartenstra\ss e 7, 64289 Darmstadt, Germany\\
Mikhail Lifshits, St.\ Petersburg State University, Russian Federation, 199034,
St.Petersburg, Universitetskaya emb. 7-9.

\end{document}